\documentclass[10pt]{amsart}
\usepackage[margin=1in]{geometry}
\usepackage[dvipsnames,usenames]{color}
\usepackage{hyperref}
\usepackage{graphicx}
\usepackage{epsfig}
\usepackage{comment}
\usepackage[latin1]{inputenc}
\usepackage{amsmath}
\usepackage{amsfonts}
\usepackage{amssymb}
\usepackage{amsthm}
\usepackage{amscd}
\usepackage{verbatim}
\usepackage{subfigure}
\usepackage{pinlabel}
\usepackage{stmaryrd}
\usepackage{enumerate, enumitem}

\usepackage[textsize=scriptsize]{todonotes}

\usepackage{tikz}
\usetikzlibrary{cd}
\usetikzlibrary{arrows}
\usetikzlibrary{decorations.pathreplacing}
\usepackage{verbatim}

\theoremstyle{definition}

\theoremstyle{remark}

\title{The unknotting number of $11n102$ is 2}

\author[Tye Lidman]{Tye Lidman}
\thanks{TL was partially supported by NSF grant DMS-2506277 and a Simons Travel Award.  He thanks Chuck Livingston and Allison Moore for their immensely valuable work on KnotInfo.  No AI was harmed during the writing of this note.}
\address{Department of Mathematics, North Carolina State University, Raleigh, NC 27607}
\email{tlid@math.ncsu.edu}

\numberwithin{equation}{section}

\newcommand\HFhat{\widehat{HF}}
\newcommand\mfs{\mathfrak{s}}

\begin{document}
\maketitle


We prove that the unknotting number of the knot $11n102$ is 2.  It is known, see e.g. KnotInfo \cite{KnotInfo}, that the unknotting number is 1 or 2, so we will assume for a contradiction it is 1.  
By the Montesinos trick \cite{Montesinos}, the branched double cover of $11n102$ is half-integral surgery on some knot $J$ in $S^3$.  Up to mirroring, the knot $11n102$ is the determinant 3 Montesinos knot $M(-2/3, 1/3, 2/7)$  \cite{KnotInfo}; with this orientation, the branched double cover $Y$ is the Seifert fibered space $S^2(0; -2/3, 1/3, 2/7)$ \cite{Montesinos}, which has $H_1 = \mathbb{Z}/3$.  We can alternatively write this manifold as $S^2(1; 1/3, 1/3, 2/7)$, where our conventions mean that in the keychain link surgery diagram, the central circle has surgery coefficient $1$ and the three other components have surgery coefficients $3, 3, 7/2$.  
From this, it is easy to check that the associated plumbing bounding $Y$ is positive-definite.  Consequently, the Heegaard Floer homology of $Y$ is readily computed using the algorithms of N\'emethi \cite{Nemethi} or Ozsv\'ath-Szab\'o \cite{OSzPlumbed}: there are two spin$^c$ structures with $HF_{red} = 0$ and $d = 1/6$, while the third has $HF_{red} = \mathbb{F}^2$ and $d = -1/2$.  In principle, $Y$ could be either $+3/2$-surgery or $-3/2$-surgery on $J$.  The mod 2 $d$-invariants of $Y$ agree with those of $+3/2$-surgery on the unknot and not $-3/2$ surgery \cite[Proposition 4.8]{OSzAbsGraded}, so the Ni-Wu formula \cite[Proposition 1.6]{NiWu} implies that $Y$ must be $+3/2$ surgery on $J$.  (The sign of the surgery coefficient can also be deduced by computing the linking form of $Y$.)   In fact, the $d$-invariants of $Y$ agree with $L(3,2)$, so in Ni-Wu's notation we have $V_0(J) = 0$ and hence $V_s(J) = H_{-s}(J) = 0$ for all $s \geq 0$ \cite[Theorem 2.5, Lemma 2.7]{NiWu}.  Because the surgery is positive, we have from \cite{NiWu} and Gainullin's work \cite[Corollary 14, Proposition 15]{Gainullin} that $HF_{red}(Y,\mfs_i) \cong \mathbb{A}^{red}_{i,3/2}$; here, the notation means there is an ordering of the three spin$^c$ structures on $Y$, denoted $\mfs_0, \mfs_1, \mfs_2$, and a sequence of vector spaces $Q_s$ for $s \in \mathbb{Z}$ so that $\mathbb{A}^{red}_{i,3/2} = \bigoplus_{k \in \mathbb{Z}}Q_{\lfloor \frac{i+3k}{2} \rfloor}$ for each $i$.  (Only finitely many $Q_s$ are non-zero.)
For any $i$ and $k$, note that $\lfloor \frac{i+3k}{2} \rfloor$ equals either $\lfloor \frac{i+1 + 3k}{2} \rfloor$ or $\lfloor \frac{i - 1+3k}{2} \rfloor$.  Therefore, each $Q_s$ shows up in two of the three $\mathbb{A}^{red}_{i,3/2}$.  This means it is impossible for exactly one of the $HF_{red}(Y,\mfs_i)$ to be non-zero.  This gives the desired contradiction.

\bibliographystyle{alpha}
\bibliography{references}

\end{document}